\documentclass[11pt]{article} 
\usepackage{amsmath,amsfonts,amssymb,mathtools}
\usepackage{graphicx}
\usepackage{hyperref}
\usepackage[italicdiff]{physics}
\usepackage{amsrefs}

\providecommand{\beq}{\begin{equation}}
\providecommand{\eeq}{\end{equation}}

\newcommand{\pint}{\mathbb Z^{+}}
\newcommand{\npint}{\mathbb Z^-_{0}}

\newcommand{\K}{K(z)}

\usepackage{authblk}

\date{}                     

\title{A new entire factorial function}
\author[ ]{Matthew D. Klimek}
\affil[ ]{Laboratory for Elementary Particle Physics, Cornell University, Ithaca, NY, 14853, USA\protect\\
Department of Physics, Korea University, Seoul 02841, Republic of Korea}

\begin{document}
\maketitle
\begin{abstract}
We introduce a new factorial function which agrees with the usual Euler gamma function at both the positive integers and at all half-integers, but which is also entire. We describe the basic features of this function.
\end{abstract}

%
The canonical extension of the factorial to the complex plane is given by Euler's gamma function $\Gamma(z)$.
It is the only such extension that satisfies the recurrence relation $z\Gamma(z)=\Gamma(z+1)$ and is logarithmically convex for all positive real $z$ \cite{mollerup,artin}.
(For an alternative function theoretic characterization of $\Gamma(z)$, see \cite{wielandt}.)
As a consequence of the recurrence relation, $\Gamma(z)$ is a meromorphic function when continued onto the complex plane. 
It has poles at all non-positive integer values of $z\in\npint$.
However, dropping the two conditions above, there are infinitely many other extensions of the factorial that could be constructed.

Perhaps the second best known factorial function after Euler's is that of Hadamard \cite{hadamard} which can be expressed as
\beq\label{hadamard}
	H(z) = \Gamma(z) \qty[ 1+\frac{\sin \pi z}{2\pi} \qty( \psi\qty(\frac{z}{2}) - \psi\qty(\frac{z+1}{2}))]
		\equiv \Gamma(z) \qty[ 1+Q(z) ]
\eeq
where $\psi(z)$ is the digamma function 
\beq\label{digamma}
	\psi(z) = \dv{z} \log\Gamma(z) = \frac{\Gamma'(z)}{\Gamma(z)} .
\eeq
The second term in brackets is given the name $Q(z)$ for later convenience.
We see that the Hadamard gamma is a certain multiplicative modification to the Euler gamma.
However, the poles are no longer present. 
It is an entire function.

In order to gain a clear picture of what the Hadamard gamma is doing, let us first note a property of the digamma function $\psi(z)$.
Since $\Gamma(z)$ has a first order pole at all $z\in\npint$, $\psi(z)$ will have first order poles at the same locations and the residues of these poles will all be $-1$.
In \eqref{hadamard} the digammas are multiplied by $\sin \pi z$, which has a first order zero at every $z\in\mathbb Z$. 
The digammas are regular at all positive integers $z\in\pint$, so $Q(z)$ vanishes for $z\in\pint$.
This is essential, because at positive integers, the Euler gamma returns the usual factorial values, and this behavior must not be spoiled.
For $z\in\npint$, one or the other of the digammas has a pole of residue $-1$, which cancels against the zero of the sine so that $Q(z)= -1$.
This in turn cancels against the first term in the square brackets, so that that factor has zeros for every $z\in\npint$.
Of course, this is precisely where the Euler gamma has poles.
So we see that the new factor in the Hadamard gamma has precisely the structure required to eliminate the poles in the Euler gamma.

The Hadamard gamma has a fairly complicated form, involving two different digamma functions.
The reason for this can be traced to the requirement that $Q(z)=-1$ for $z\in\npint$.
The poles of the digamma always have residue $-1$, 
but the sine approaches zero from different directions at each successive integer, so that the product of the sine with a single digamma would alternate between $\pm 1$.
$Q(z)$ therefore includes two digammas with opposite sign with arguments such that one has a pole for even non-positive integers and the other for odd non-positive integers.

\begin{figure}
	\begin{center}
		\includegraphics[width=0.5 \linewidth]{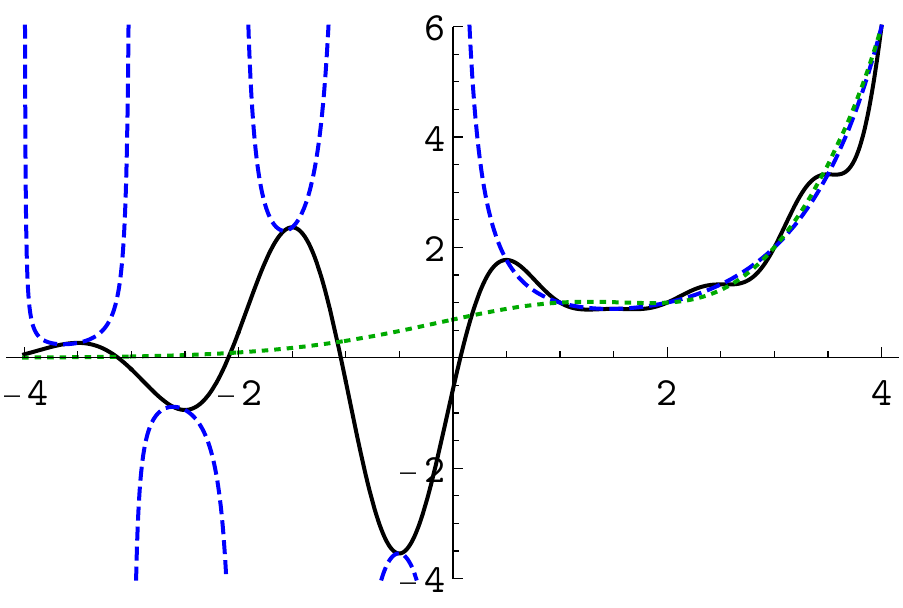}
		\includegraphics[width=0.4 \linewidth]{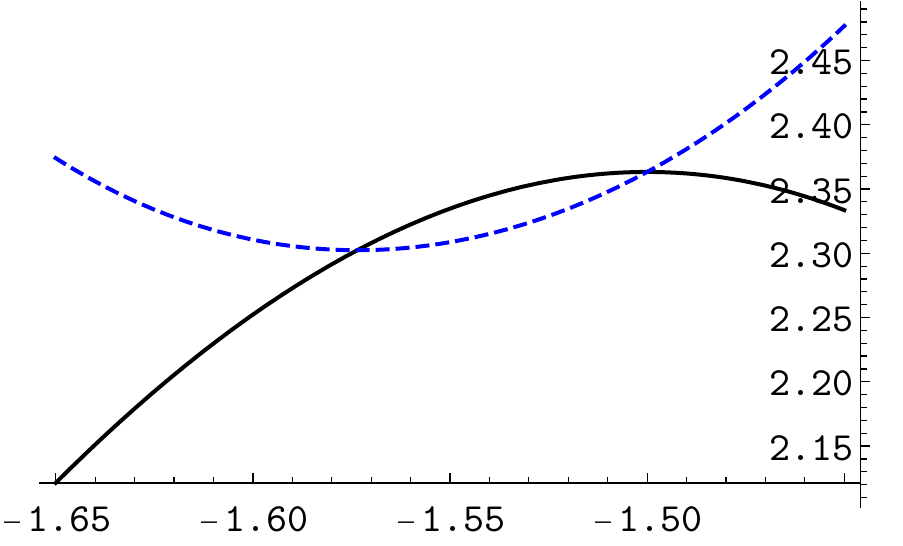}
	\end{center}
	\caption{\emph{Left:} $\K$ is plotted (black, solid) along with the Euler gamma $\Gamma(x)$ (blue, dashed) and the Hadamard gamma $H(x)$ (green, dotted).
	\emph{Right:} An example of the phenomenon that $\K$ and the Euler gamma are equal wherever one of them has an extremum.}
\label{fig:plots} 
\end{figure}

\newcommand{\sintwopi}{\frac{\sin 2\pi z}{2\pi}}
This observation makes it immediately obvious how to construct a simpler factorial function.
Doubling the frequency of the sine, it will always approach zero from the same direction at the integers and only one digamma is needed.
We therefore define a factorial function $K(z)$ with only a single digamma as
\beq\label{kfunction}
	K(z) = \Gamma(z) \qty( 1+ \sintwopi \psi(z) ) .
\eeq
Alternatively, $\K$ may be written as
\beq
	\K = \Gamma(z)+\sintwopi\,\Gamma'(z).
\eeq
In the second term, $\Gamma'(z)$ has second order poles for all $z\in\npint$. These are reduced to first order poles by the sine, so that the second term is subtracting the poles from the original gamma function.
This function along with the Euler and Hadamard gammas are plotted in the left panel of Figure~\ref{fig:plots}.

Like the Hadamard gamma, $\K$ is an entire function.
It also has a number of other interesting properties.
Due to the doubling of the frequency of the sine, $\K$ agrees with the Euler gamma at all half-integer values as well.
As half-integer values of the factorial arise naturally in many contexts, this may be seen as an attractive feature.
It also agrees with the Euler gamma at all of the Euler gamma's extrema, or equivalently, at the zeros of the digamma function.
The derivative may be written as
\beq
	K'(z) = \Gamma(z)\qty[ \qty(1+\cos 2\pi z)\psi(z) + \frac{\sin 2\pi z}{2\pi}\qty(\psi^2(z)+\psi'(z))].
\eeq
From this we see the remarkable fact that the derivative vanishes at every half-integer value of $z$, where $\K$ also equals the Euler gamma.
We conclude that $\K=\Gamma(z)$ anywhere that either of them has an extremum.
This is illustrated near the point $z=-3/2$ in the right panel of Figure~\ref{fig:plots}.
Meanwhile, for $z\in\pint$ where $\K=(z-1)!$ we have $K'(z)=2\Gamma'(z)$.

The recurrence relation satisfied by $\K$ is the same as that of $\Gamma(z)$ plus an extra term which vanishes at all integers and half-integers
\beq
	K(z+1) = z\K +\Gamma(z)\,\sintwopi .
\eeq
This can also be written as
\beq
	K(z+1) = z\K +\frac{\cos\pi z}{\Gamma(1-z)}
\eeq
by applying the reflection identity $\Gamma(z)\Gamma(1-z)=\pi/\sin\pi z$.
This is similar to the relation 
\beq\label{Hrecurrence}
	H(z+1) = zH(z) +\frac{1}{\Gamma(1-z)}
\eeq
for the Hadamard gamma but with an oscillating factor in the second term.
In any case, for both $\K$ and the Hadamard gamma, since $\Gamma(1-z)$ grows rapidly for $x\to -\infty$, these relations revert to the original recurrence relation of the Euler gamma asymptotically, but no poles appear due to the modification of the recurrence relation near the origin. 

At $z=0$, $\K$ takes the finite value $K(0)=-\gamma$ where $\gamma$ is the Euler-Mascheroni constant.
Going to positive values of $z$, $\K$ has a zero at $z=0.06958\dots$ before rising to $K(1/2)=\Gamma(1/2)=\sqrt\pi$.
From here, $\K$ continues to increase to match the values of $\Gamma(z)$ at (half-)integer values, while oscillating around $\Gamma(z)$.
For small values of $z$, the amplitude of these oscillations is small so that $\K$ has no additional zeros.
However, $\psi(z)$ has the asymptotic expansion $\psi(z)\sim\log z +\order{z^{-1}}$ so that starting at values of $z$ of order $\exp( 2\pi)$ additional pairs of zeros appear.
This violent oscillation at large values of the argument may seem unattractive, but it is a consequence of the simpler form of $\K$ and the fact that it also reproduces half-integer factorials.

Finally, we note that one can construct an class of entire functions that match the Euler gamma at integers and half-integers by adding to $\K$ a term of the form
\beq
	e^{g(z)} \Gamma(z) \sin 2\pi z
\eeq
where $g(z)$ is an arbitrary entire function.
This follows from the fact that $\Gamma(z)\sin 2\pi z$ is an entire function with zeros at all  half-integers and positive integers, so adding it to $\K$ does not change its agreement with the factorial at those points, and that any two entire functions $f_1(z)$ and $f_2(z)$ with the same zeros can be related to each other as $f_1(z) = e^{g(z)} f_2(z)$ for a suitable entire function $g(z)$.

\section*{Acknowledgements}
The author acknowledges the support of the U.S. National Science Foundation, through grant PHY-1719877, the Samsung Science \& Technology Foundation, under Project Number SSTF-BA1601-07, and a Korea University Grant.


\begin{bibdiv}\begin{biblist}

\bib{mollerup}{book}{
title={L{\ae}rebog i Kompleks Analyse },
author={Mollerup, J.},
author={Bohr, H.},
volume={3},
date={1922},
publisher={Copenhagen}
}

\bib{artin}{book}{
title={Einf\"uhrung in die Theorie der Gammafunktion},
author={Artin, E.},
date={1931},
publisher={Teubner},
translation={
title={The Gamma Function},
translator={Butler, M.},
date={1964},
publisher={Holt, Rinehart and Winston}
}
}

\bib{wielandt}{article}{
title={Wielandt's Theorem About the $\Gamma$-Function },
author={Remmert, R.},
journal={American Mathematical Monthly },
volume={103},
number={3},
date={1996},
pages={214--220}
}

\bib{hadamard}{article}{
title={ Sur L'Expression Du Produit $1\cdot2\cdot3\cdots(n-1)$ Par Une Fonction Enti\`ere},
author={Hadamard, M. J.},
booktitle={{\OE}uvres de Jacques Hadamard },
publisher={Centre National de la Recherche Scientifiques},
date={1968},
}

\end{biblist}\end{bibdiv}
\end{document}